\documentclass[reqno]{amsart}
\numberwithin{equation}{section}

\usepackage{amsmath}
\usepackage{amsthm}
\usepackage{amssymb}
\usepackage{graphics}

\usepackage{mathrsfs} 

\usepackage{comment} 

\newcommand{\R}{\mathbb{R}}

\usepackage{bm}
\usepackage{graphicx}
\usepackage{amscd}
\usepackage{pstricks}
\usepackage{pst-plot}
\usepackage{multido}
\usepackage{pst-coil}
\usepackage{colortab}

\usepackage{comment} 

\usepackage{tikz} 

\usepackage{float} 

\theoremstyle{plain}
\newtheorem{theorem}{\bf Theorem}[section]

\theoremstyle{remark}
\newtheorem{remark}[theorem]{\bf Remark}

\usepackage{float}
\begin{document}


\title[Addition theorem for 3-D Navier-Lamé system and applications]{A new addition theorem for the 3-D Navier-Lamé system and its application to the method of fundamental solutions}
\author{J.A. Barceló, C. Castro, A. Ruiz and M.C. Vilela}
\newcommand{\Addresses}{{
  \bigskip
  \footnotesize

  J.A.~Barceló, 
  \textsc{Universidad Politécnica de Madrid, Spain}\par\nopagebreak
  \textit{E-mail address}: \texttt{juanantonio.barcelo@upm.es}

  C.~Castro, 
  \textsc{Universidad Politécnica de Madrid, Spain}\par\nopagebreak
  \textit{E-mail address}: \texttt{carlos.castro@upm.es}

  \medskip

  A.~Ruiz, 
  \textsc{Universidad Autónoma de Madrid, Spain}\par\nopagebreak
  \textit{E-mail address}: \texttt{alberto.ruiz@uam.es}

  \medskip

  M.C.~Vilela, 
  \textsc{Universidad Politécnica de Madrid, Spain}\par\nopagebreak
  \textit{E-mail address}: \texttt{maricruz.vilela@upm.es}

}}

\thanks{The first, the second and the fourth authors were supported by the Spanish Grants PID2021-124195NB-C31 and PID2024-158664NB-C21. The third one is supported by the Spanish Grants PID2021-124195NB-C32 and PID2024-158664NB-C22.}

\begin{abstract}
    We derive a new addition formula for the fundamental solution of the Navier–Lamé system in three dimensions that satisfies the Kupradze radiation conditions. This result yields an expansion of the fundamental solution involving only Bessel functions and scalar spherical harmonics. Such a representation is particularly advantageous for collocation-based numerical methods relying on fundamental solutions, including the boundary element method and the method of fundamental solutions. For the latter, we illustrate its efficiency in approximating solutions to the Navier–Lamé system in exterior domains.
\end{abstract}

\maketitle

\section{Introduction}
In the present paper, we prove an addition formula for the fundamental solution of the Navier-Lamé system in three dimensions with Kupradze radiation conditions at infinity. This is an extension of the analogous formula obtained for the two dimensional system in \cite{BRVW2024}. Addition theorems are an important tool in obtaining a priori estimates for solutions of partial differential equations (see \cite{HartmanWilcox} and \cite[Theorem 1]{BRV1997} for the Helmholtz equation or  \cite{ColtonKress} for the Maxwell system), but also in numerical methods that require explicit formulas for the fundamental solutions, as trigonometric collocation methods (see \cite{BarceloCastro2022}), boundary elements or the so-called method of fundamental solutions (MFS). We refer to \cite{ChengHong2020} for an exhaustive review of this last method in different situations.  

In contrast to previously known expressions of the fundamental solution, our addition formula does not involve derivatives of spherical harmonics or Bessel functions (see \cite{Kupradze1979} and \cite{HsiaoWendland}). This is important for numerical methods based on collocation, since they require a large number of evaluations at specific points. To illustrate this, the MFS is considered to solve the Navier-Lamé system in an exterior domain with boundary data. We show that the method provides accurate solutions with a relatively small number of terms in the asymptotic expansion.  We give numerical experiments for both the dimension $d=2$ and $3$.

The application of the addition theorem to obtain a priori estimates for the solution of the Navier-Lamé system is described in \cite[Theorem 1]{BRVW2024}, where an analogous result is established in two dimensions.

The rest of this paper is organized as follows: in Section 2 we state the addition theorem for the fundamental solution of the Naver-Lamé system with Kupradze radiation conditions. This is an explicit series that only requires evaluation of the Hankel and Bessel functions. In Section 3 we show how we can use this to solve the boundary value problem in a exterior domain. In Section 4 we give a proof of the result stated in Section 2.  In Section 5 we introduce several technical identities required for this proof. Finally, in Section 6 we discuss the main conclusions.
 
\section{The addition theorem for the fundamental solution of the Navier-Lamé system}

Consider the Navier-Lamé system in the three dimensional case 
\begin{equation}
\Delta^{\ast}\mathbf{u}(x)+\omega^{2}\mathbf{u}(x) = \mathbf{f}(x), \qquad\qquad \omega>0,\ x\in\mathbb{R}^{3}, 
\label{ecuacion}
\end{equation}
where $\mathbf{u}$ is a vector-valued function from
$\mathbb{R}^{3}$ to $\mathbb{R}^{3}$ and,
\begin{equation}
\label{operador}\Delta^{\ast}\mathbf{u}(x) = \mu\Delta\mathbf{u}(x)+(\lambda+\mu)\nabla div\,\mathbf{u}(x),
\end{equation}
with $\Delta\mathbf{u}$ defined component by component and, $\mu$ and $\lambda$ are the Lamé constants.

We will assume that $\mu>0$ and
$2\mu+\lambda>0$ so that the operator $\Delta^{\ast}$ is strongly
elliptic and, we will denote by $k_p$ and $k_s$  respectively the speed of propagation of  longitudinal and  transversal waves, where
\begin{equation}
\label{kpks} \displaystyle k_p^2=\frac{\omega^2}{(2\mu+\lambda)}
\qquad \texttt{\rm{and}} \qquad k_s^2=\frac{\omega^2}{\mu}.
\end{equation}

The main result of this paper is the following expansion for the fundamental solution of \eqref{ecuacion} satisfying the Kupradze radiation conditions (see \cite{Kupradze}).


\begin{theorem}
    \label{th_adition_theorem}
    Let $Y_n^m$, $m=-n,\ldots,n$, $n=0,1,\ldots$, be the set of orthonormal spherical harmonics given in \eqref{gamma_n^m}. 
    Denote the spherical coordinates of $x$ and $y$ by $(r,\theta_x,\varphi_x)$ and $(t,\theta_y,\varphi_y)$ respectively, where \( r = |x| \), \( t = |y| \), \( \theta_x, \theta_y \in [0,\pi] \) are the polar angles, and \( \phi_x, \phi_y \in [0,2\pi) \) are the azimuthal angles.
	Then for $|x|>|y|$, the fundamental solution of \eqref{ecuacion} satisfying the Kupradze radiation conditions is a matrix that can be written as
 \begin{equation}
 \label{additionLame}
    \mathrm{\Phi}(x,y)=\mathrm{\Psi}(x,y)+\mathrm{\Phi}_-(x,y)+\mathrm{\Phi}_0(x,y)+\mathrm{\Phi}_+(x,y),
\end{equation}
with
    \begin{align}
	\mathrm{\Psi}(x,y)&=
    \frac{i}{\omega^2}\sum_{n=0}^\infty
    (H_{n,n}^{k_s}(r,t)\mathrm{S}
    +
    \frac{1}{2}H_{n,n}^{+,k_p,k_s}(r,t) \mathrm{P})
    \sum_{m=-n}^{n}
    Y_{n}^{m}(\theta_x,\varphi_x)\overline{Y_n^m}(\theta_y,\varphi_y),
\label{Phi1_3D}
    \\
   \mathrm{\Phi}_-(x,y)&=
    \frac{1}{\omega^2}\sum_{n=0}^\infty
    H_{n-2,n}^{-,k_p,k_s}(r,t)
    \sum_{m=-n}^{n}
    S_{-,n}^m(\theta_x,\varphi_x)\overline{Y_n^m}(\theta_y,\varphi_y),
\label{Phi2_3D}
    \\
    \mathrm{\Phi}_0(x,y)&=
    \frac{1}{\omega^2}\sum_{n=0}^\infty
    H_{n,n}^{-,k_p,k_s}(r,t)
    \sum_{m=-n}^{n}
    S_{0,n}^m(\theta_x,\varphi_x)\overline{Y_n^m}(\theta_y,\varphi_y),
\label{Phi3_3D}
    \\
    \mathrm{\Phi}_+(x,y)&=
    \frac{1}{\omega^2}\sum_{n=0}^\infty
    H_{n+2,n}^{-,k_p,k_s}(r,t)
    \sum_{m=-n}^{n}
    S_{+,n}^m(\theta_x,\varphi_x)\overline{Y_n^m}(\theta_y,\varphi_y),
\label{Phi4_3D}
    \end{align}
where for $n_1,n_2\ge 0$
\begin{align}
\label{h0}
H_{n_1,n_2}^{k_s}(r,t)&:=
k_s^3h_{n_1}^{(1)}(k_sr)j_{n_2}(k_st),
\\
\label{h+}
H_{n_1,n_2}^{+,k_p,k_s}(r,t)&:=
k_p^3h_{n_1}^{(1)}(k_pr)j_{n_2}(k_pt)
+k_s^3h_{n_1}^{(1)}(k_sr)j_{n_2}(k_st),
\\
\label{h-}
H_{n_1,n_2}^{-,k_p,k_s}(r,t)&:=
k_p^3h_{n_1}^{(1)}(k_pr)j_{n_2}(k_pt)
-k_s^3h_{n_1}^{(1)}(k_sr)j_{n_2}(k_st),
\end{align}
with $j_n$ and $h_n^{(1)}$ denoting the spherical Bessel functions of the first and third kind respectively,
\begin{equation*}
    \mathrm{S}=\left(
    \begin{array}{ccc}
       0&0&0  \\
       0&0&0  \\
       0&0&1 
    \end{array}
    \right),
    \qquad
    \mathrm{P}=\left(
    \begin{array}{ccc}
       1&0&0  \\
       0&1&0  \\
       0&0&0 
    \end{array}
    \right),
\end{equation*}
\begin{align}
    S_{-,n}^m(\theta_x,\varphi_x)
    =&
    \mathrm{A}_{-,n}^{m}
    Y_{n-2}^{m-2}(\theta_x,\varphi_x)
    +
    \mathrm{B}_{-,n}^{m}
    Y_{n-2}^{m-1}(\theta_x,\varphi_x)
    +
    \mathrm{C}_{-,n}^{m}
    Y_{n-2}^{m}(\theta_x,\varphi_x)
    \nonumber
    \\
    &+
    \mathrm{D}_{-,n}^{m}
    Y_{n-2}^{m+1}(\theta_x,\varphi_x)
    +
    \mathrm{E}_{-,n}^{m}
    Y_{n-2}^{m+2}(\theta_x,\varphi_x),   
\label{S-2}
\end{align}
\begin{align}
    S_{0,n}^m(\theta_x,\varphi_x)
    =&
    \mathrm{A}_{0,n}^{m}
    Y_{n}^{m-2}(\theta_x,\varphi_x)
    +
    \mathrm{B}_{0,n}^{m}
    Y_{n}^{m-1}(\theta_x,\varphi_x)
    +
    \mathrm{C}_{0,n}^{m}
    Y_{n}^{m}(\theta_x,\varphi_x)
    \nonumber
    \\
    &+
    \mathrm{D}_{0,n}^{m}
    Y_{n}^{m+1}(\theta_x,\varphi_x)
    +
    \mathrm{E}_{0,n}^{m}
    Y_{n}^{m+2}(\theta_x,\varphi_x),   
\label{S0}
\end{align}
\begin{align}
    S_{+,n}^m(\theta_x,\varphi_x)
    =&
    \mathrm{A}_{+,n}^{m}
    Y_{n+2}^{m-2}(\theta_x,\varphi_x)
    +
    \mathrm{B}_{+,n}^{m}
    Y_{n+2}^{m-1}(\theta_x,\varphi_x)
    +
    \mathrm{C}_{+,n}^{m}
    Y_{n+2}^{m}(\theta_x,\varphi_x)
    \nonumber
    \\
    &+
    \mathrm{D}_{+,n}^{m}
    Y_{n+2}^{m+1}(\theta_x,\varphi_x)
    +
    \mathrm{E}_{+,n}^{m}
    Y_{n+2}^{m+2}(\theta_x,\varphi_x),   
\label{S-}
\end{align}
with the matrix coefficientes $\mathrm{A}_{\#,n}^{m},\mathrm{B}_{\#,n}^{m},\mathrm{C}_{\#,n}^{m},\mathrm{D}_{\#,n}^{m}$ and $\mathrm{E}_{\#,n}^{m}$ for $\#=-,+,0,$ given in the following table:
\begin{table}[H]
    \centering
    \begin{tabular}{c|c|c|c|}
        \# &$-$&$0$&$+$
        \\[0.5ex]\hline
        $\mathrm{A}_{\#,n}^{m}$ 
        & $\frac{1}{4}b_{n-2}^{m-2}c_n^m\mathrm{A}$
        & $\frac{1}{4}(b_{n}^{m-2}a_n^m+d_{n}^{m-2}c_n^m)\mathrm{A}$
        & $\frac{1}{4}d_{n+2}^{m-2}a_n^m\mathrm{A}$
        \\[0.5ex]\hline
        $\mathrm{B}_{\#,n}^{m}$
        & $\frac{1}{2}b_{n-2}^{m-1}f_n^m\mathrm{B}$
        & $\frac{1}{2}(d_{n}^{m-1}f_n^m-b_{n}^{m-1}e_n^m)\mathrm{B}$
        & $-\frac{1}{2}d_{n+2}^{m-1}e_n^m\mathrm{B}$
        \\[0.5ex]\hline
        $\mathrm{C}_{\#,n}^{m}$
        & $\frac{1}{2}a_{n-2}^{m}c_n^m\mathrm{C}$
        & $-\frac{1}{2}((e_n^m)^2+(f_n^m)^2)\mathrm{C}$
        & $\frac{1}{2}c_{n+2}^{m}a_n^m\mathrm{C}$
        \\[0.5ex]\hline
        $\mathrm{D}_{\#,n}^{m}$
        & $\frac{1}{2}a_{n-2}^{m+1}f_n^m\mathrm{D}$
        & $\frac{1}{2}(c_{n}^{m+1}f_n^m-a_{n}^{m+1}e_n^m)\mathrm{D}$
        & $-\frac{1}{2}c_{n+2}^{m+1}e_n^m\mathrm{D}$
        \\[0.5ex]\hline
        $\mathrm{E}_{\#,n}^{m}$
        & $\frac{1}{4}a_{n-2}^{m+2}d_n^m\mathrm{E}$
        & $\frac{1}{4}(a_{n}^{m+2}b_n^m+c_{n}^{m+2}d_n^m)\mathrm{E}$
        & $\frac{1}{4}c_{n+2}^{m+2}b_n^m\mathrm{E}$
        \\[0.5ex]\hline
    \end{tabular}
    \caption{Matrix coefficients}
    \label{tabla_phi}
\end{table}
where $a_{n,m},b_{n,m},c_{n,m},d_{n,m}$ and $e_{n,m}$ are functions of discrete variables $n$ and $m$, with $n\ge0,$ and $-n\le m\le n$, uniformly bounded and given
in \eqref{a}, \eqref{d} and  \eqref{e}, and $\mathrm{A},\mathrm{B}, \mathrm{C}, \mathrm{D}$ and $\mathrm{E}$ are the following matrices of complex numbers:

\begin{equation*}
    \mathrm{A}=\left(
    \begin{array}{ccc}
       -i&1&0  \\
       1&i&0  \\
       0&0&0 
    \end{array}
    \right),\ 
    \mathrm{B}=\left(
    \begin{array}{ccc}
       0&0&i  \\
       0&0&-1  \\
       i&-1&0 
    \end{array}
    \right),\ 
    \mathrm{C}=\left(
    \begin{array}{ccc}
       i&0&0  \\
       0&i&0  \\
       0&0&-2i 
    \end{array}
    \right),
\end{equation*}
    \begin{equation*}
    \mathrm{D}=\left(
    \begin{array}{ccc}
       0&0&-i  \\
       0&0&-1  \\
       -i&-1&0 
    \end{array}
    \right),\ 
    \mathrm{E}=\left(
    \begin{array}{ccc}
       -i&-1&0  \\
       -1&i&0  \\
       0&0&0 
    \end{array}
    \right).
    \end{equation*}

The series and its term by term first derivatives with respect to $|x|$ and $|y|$ are absolutely and uniformly convergent on compact subsets of $|x| > |y|$.
\end{theorem}

 \begin{remark}
     \label{venta}
     We would like to note that the functions given in \eqref{h-} reflect compensations between the longitudinal and transversal waves in the fundamental solution. This fact is crucial to obtain apriori estimates for solutions of the system \eqref{ecuacion} (see \cite[Lemmas 4.3 and 4.4]{BRVW2024}).
 \end{remark}

 \section{The MFS for an exterior domain}

 In this section we give an application of the main result in Theorem \ref{th_adition_theorem}. In particular we apply this result to approximate the elasticity system in an exterior domain using the method of fundamental solutions (MFS). The MFS is a boundary collocation method first proposed as a computational technique by Mathon and Johnston \cite{MathonJohnston1977}, building on the theoretical framework developed by Kupradze and Aleksidze \cite{KupradzeAleksidze1964}. It is one of the simplest numerical methods based on explicit representations of fundamental solutions, and has been shown to be particularly effective for exterior and unbounded domain problems \cite{FairweatherKarageorghis1998}. 
Unlike previously known expressions of the fundamental solution of the
elasticity operator, which are given in terms of derivatives of Hankel and
Bessel functions \cite{Kupradze1979}, the addition formula derived in
Theorem~\ref{th_adition_theorem} requires only evaluations of these
functions, which is particularly advantageous for collocation methods such
as the MFS.
 
 We consider as domain $\Omega$ the exterior of the cube of side length 4 centered at the origin. 
 
 We are interested in approximating the solution of the problem  
\begin{equation} \label{eq_BVP}
    \left\{ 
    \begin{array}{ll}
        \Delta^{\ast}\mathbf{u} +\omega^{2}\mathbf{u}=0  & \mbox{ in $\Omega$} \\
         \mathbf{u}=\mathbf{g} & \mbox{ on $\partial \Omega$},
    \end{array}
    \right.
\end{equation}
for a given continuous function $\mathbf{g}$ defined on the boundary $\partial \Omega$. The idea of MFS is to approximate $\mathbf{u}$ as a linear combination of fundamental solutions with singularities outside $\Omega$, i.e.
\begin{equation}
  \mathbf{u}(x)\sim \mathbf{u}^N(x)=\sum_{k=1}^N  \Phi(x,y_k) \bm{\alpha}_k,  
\end{equation}
where $\{y_k\}_{k=1}^N\subset \R^3 \backslash \Omega$ are some a priori chosen basis points. The unknown vector coefficients $\bm{\alpha}_k$  are determined by collocation on the boundary. This requires choosing some boundary points $\{ x_j \}_{j=1}^M\subset \partial \Omega$ where the boundary condition is assumed to be satisfied, i.e.
\begin{equation} \label{eq_colo}
  \mathbf{u}^N(x_j)=\mathbf{g}(x_j), \quad j=1,...,M.  
\end{equation}

Equations \eqref{eq_colo} are then reduced to the linear system 
\begin{equation} \label{eq_SL}
    \mathbf{M}\bm{\alpha}^N= \mathbf{g}^M,
\end{equation}
where $\mathbf{M}\in \mathcal{M}_{3M\times 3N}$, is given by 
\begin{equation} \label{eq_mat}
   \mathbf{M}=(\Phi(x_j,y_k))_{jk}, \quad  j=1,...,M, \quad k=1,...,N,
\end{equation}
$\bm{\alpha}^N=(\bm{\alpha}_1,...,\bm{\alpha}_N)^T \in \R^{3N}$ and $\mathbf{g}^M=(\mathbf{g}(x_1),...,\mathbf{g}(x_M))^T \in \R^{3M}$.

When $N=M$ the matrix $\mathbf{M}$ is square, and the system may have a unique solution. The cases $N>M$ or $M>N$ produce underdetermined and overdetermined systems, respectively. In these cases a least squares method can be used to construct an approximation. 

The efficiency of the MFS strongly depends on the choice of the basis points $\{y_k\}_{k=1}^N\subset \R^3 \backslash \Omega$ and the collocation points $\{ x_j \}_{j=1}^M\subset \partial \Omega$ (see \cite{FairweatherKarageorghis1998, ChengHong2020}). Associated to the set of points $\{y_k\}_{k=1}^N$ we construct a basis on a subspace $X^N\subset L^2(\partial \Omega)^9$ given by 
\begin{equation} \label{eq_basis}
X^N=span \; \{ \Phi(x,y_1)|_{\partial \Omega},...,\Phi(x,y_N)|_{\partial \Omega} \} .    
\end{equation}
The density of such MFS approximation spaces in $L^2(\partial\Omega)$
has been studied in \cite{Bogomolny1985, KatsuradaOkamoto1988} for the Hemholtz equation. 
The placement of the source points can be interpreted as determining the approximation space $X^N$, and therefore plays a fundamental role in the approximation properties of the MFS, although optimal choices are, in general, not known in a rigorous sense. In special cases, however, such as when $\Omega$ is a disk, where the source points are chosen equidistantly on a concentric circle, and the boundary data are analytic, the exponential convergence of the method is known both for the Laplace and Hemholtz equations \cite{Bogomolny1985, Katsurada1990, BarnettBetcke2008}.
On the other hand, the choice of the collocation points $\{ x_j \}_{j=1}^M$ affects the condition number of the matrix $\mathbf{M}$  (see \cite{FairweatherKarageorghis1998, BarnettBetcke2008}). 

\begin{figure}[htbp]
    \centering
    \includegraphics[width=0.5\textwidth]{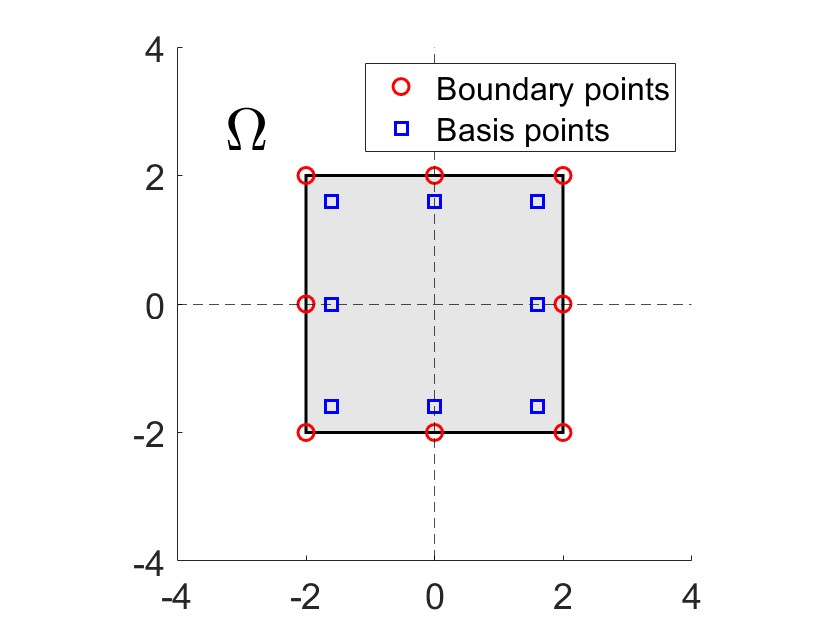}
    \caption{Distribution of boundary points $x_k$ (red circles) and basis points $y_j$ (blue squares) in the 2-d case. In this paper we consider the analogous in 3-d.}
    \label{fig1}
\end{figure}

In this paper we do not address all these issues since our main objective is to show the applicability of the addition formula stated in Theorem \ref{th_adition_theorem} for the MFS. We consider a particular example with $\{ x_j \}_{j=1}^M$ uniformly  distributed on the boundary $\partial \Omega$ and the same number of basis points $\{y_k\}_{k=1}^M$ defined by homothetically with ratio $a=0.95$
\begin{equation} \label{eq_homo}
    y_j=0.95 \; x_j,  \quad j=1,...,M=N. 
\end{equation}
In Figure \ref{fig1} we show the distribution of the points in the two-dimensional case.

\begin{figure}[htbp]
    \centering
    \includegraphics[width=0.8\textwidth]{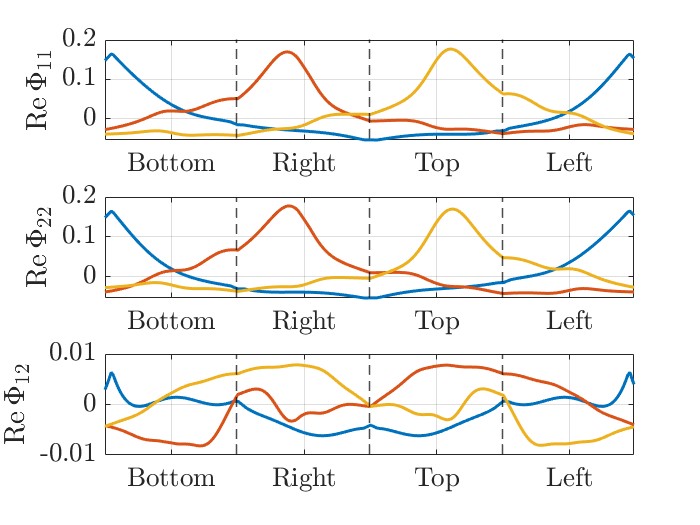}
    \includegraphics[width=0.8\textwidth]{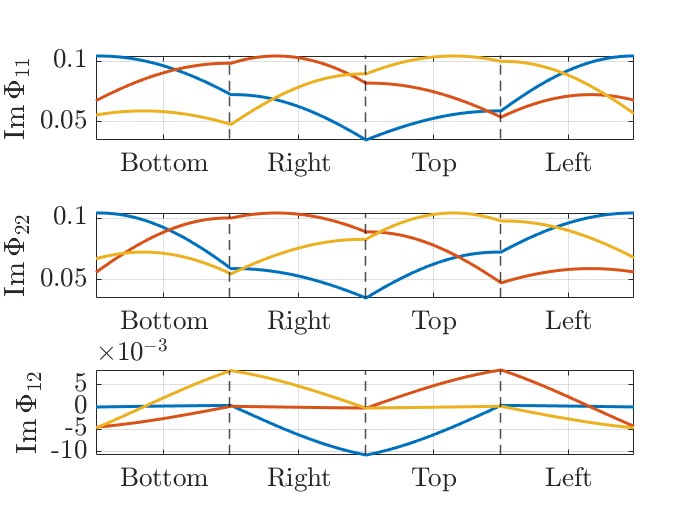}
    
    \caption{Real and imaginary parts of basis functions of the form $\Phi(x,y_j)|_{\partial \Omega} \in X^N$ corresponding to three different points $y_j$, the first one close to the bottom left corner (blue), the second one close to the right boundary (red) and the third one close to the top boundary (yellow). The bottom, right, top and left parts correspond to the four sides of the square boundary.}
    \label{fig1b}
\end{figure}

We also show in Figure~\ref{fig1b} the real and imaginary parts of the components
of three elements of the basis $X^N$ in~\eqref{eq_basis} corresponding to three different
source points $y_j$, in dimension $d=2$ with $\lambda=-1$, $\mu=2$ and $\omega=1$. Observe that the real and imaginary parts of the diagonal terms are larger near the boundary points closest
to $y_j$. This behavior is consistent with that of the Helmholtz equation, which is expected since system \eqref{operador} reduces to a vector Helmholtz equation when $\lambda+\mu=0$.

To check the approximation we consider an experiment where we take an explicit solution \textbf{u}  of the problem \eqref{eq_BVP}, so that we can compute the boundary value $\textbf{g}= \mathbf{u}|_{\partial \Omega}$ and check the approximation on a set of points in $\Omega$. A simple choice is  
    \begin{equation} \label{exact}
       \mathbf{u}(x)=\Phi(x,P)\mathbf{v}  
    \end{equation}
    with $\mathbf{v}\in \R^3$ a constant vector and $P\notin \Omega$, so that $\mathbf{u}$ is solution of \eqref{eq_BVP} with 
    \begin{equation} \label{eq_dd}
       \mathbf{g}=\Phi(x,P)\mathbf{v}|_{\partial \Omega}. 
    \end{equation}

For completeness we give experiments in dimension $d=2$ and $3$. The analogous to the addition formula stated in Theorem \ref{th_adition_theorem} for dimension $d=2$ is given in \cite{BRVW2024}. We divide the remainder of this section into two subsections, in which we analyze separately the cases $d=2$ and $d=3$.

\subsection{Dimension $d=2$}

\begin{figure}[ht]
    \centering
    \includegraphics[width=0.68\textwidth]{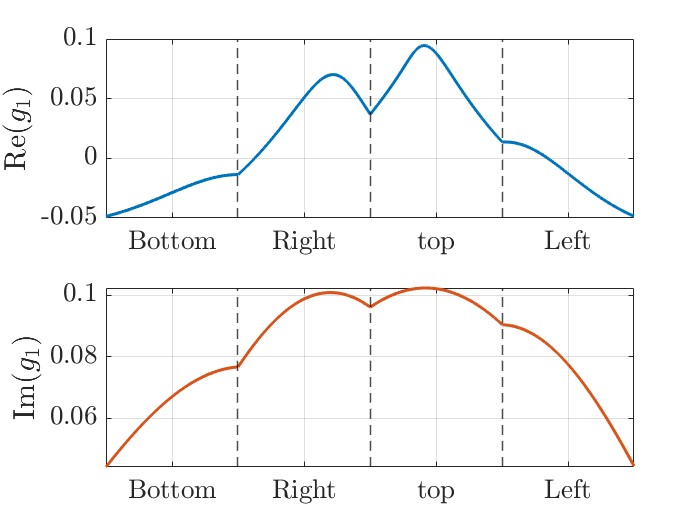}
    \caption{Real and imaginary part of the first component of the Dirichlet data $\mathbf{g}$ on the boundary square.}
    \label{fig:ambas1}
\end{figure}

We take $P$ the point with polar coordinates $(r_P,\theta_P)=(0.7,\pi/3)$ and $\mathbf{v}=(1,2)^T$. We also assume $\lambda=-1$, $\mu=2$, $\omega=1$ and approximate the exact solution in \eqref{exact} truncating the series in the $2-d$ analogous to \eqref{Phi1_3D}-\eqref{Phi4_3D} (see \cite{BRVW2024}), with a large number of terms $N_{terms}=80$ to have good accuracy. We take this as the exact solution $\mathbf{u}$ that we use both to compute the Dirichlet data $\mathbf{g}$  in \eqref{eq_dd} (see Figure \ref{fig:ambas1}) and to test the approximate solution. 

Now, to find the approximation $\mathbf{u}^N$ we follow the steps:
\begin{itemize}
    \item[Step 1.] Define $N=M$ boundary  and basis points. In our experiments, we consider a uniform mesh on the boundary for the boundary points and the corresponding homothetic points with the ratio $a=0.95$ for the basis points, as described in \eqref{eq_homo}. This configuration is presented in Figure \ref{fig1}. 

    \item[Step 2.] Construct the matrix $\mathbf{M}$ in \eqref{eq_mat}, taking the terms $N_{terms}$ to approximate the fundamental solution $\Phi$, and the vector $\mathbf{g}^N$ from the exact solution.  

    \item[Step 3.] Solve the linear system \eqref{eq_SL}. This gives the approximate solution $\mathbf{u}^N$
\end{itemize}
    
    Finally, to evaluate the error we consider the maximum of the relative  difference in a grid of points $\{ z_l\}_{l=1}^L\subset \Omega$, i.e. 
    $$
    e_\infty = \max_{l} \frac{\| \mathbf{u}(z_l)-\mathbf{u}^N(z_l) \|}{\| \mathbf{u}(z_l) \|}.
    $$
    Here we take $\{ z_l\}_{l=1}^L\subset \Omega$ those points of a uniform grid with $500$ points on the cube $[-5,5]^2$ that are in $\Omega$.

    \begin{figure}[ht]
    \centering
    \includegraphics[width=0.48\textwidth]{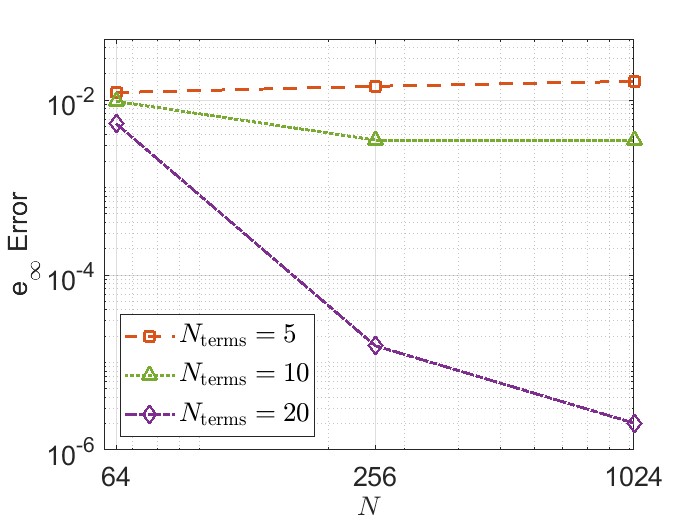}
    \hfill
    \includegraphics[width=0.48\textwidth]{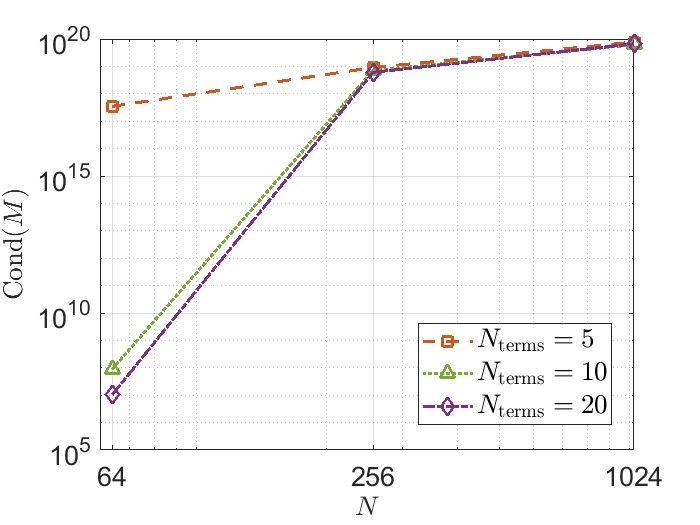}
    \caption{Dimension $d=2$. Left: $e_\infty$-error versus the number of boundary collocation points when considering different truncations of the fundamental solution up to $N_{terms}$. Right: condition number of the matrix $\mathbf{M}$.}
    \label{fig:ambas}
\end{figure}
In Figure \ref{fig:ambas} we show the $e_\infty$ error in terms of the number of collocation points in the boundary $N$ and the truncation $N_{terms}$ in the series of the fundamental solution. We observe that, for a relatively small number in $N$ and $N_{terms}$ the errors become very small. Moreover, as expected, the precision increases when both $N$ and $N_{terms}$ grow.  At the same time, the condition number of the matrix $\mathbf{M}$ grows rapidly with the number of points. This makes the linear system  \eqref{eq_SL} difficult to solve accurately when $N$ is relatively large. This is in fact one of the main drawbacks associated to the MFS. Several strategies have been proposed in the literature to mitigate this issue, such as Tikhonov regularization or SVD-based least-squares solvers; a systematic study of these techniques in this context is left for future work.

\subsection{Dimension $d=3$}
    We now take $P$ the point with spherical coordinates $(r_P,\theta_P,\varphi_P)=(0.7,1,1)$ and $\mathbf{v}=(1,2,-1)^T$. We also assume $\lambda=-1$, $\mu=2$ and $\omega=1$. As in the 2D case, we approximate the exact solution by truncating the series in \eqref{Phi1_3D}-\eqref{Phi4_3D} with $Nterms=40$, which provides sufficient accuracy to serve as a reference solution for the error computation.
    
\begin{figure}[ht]
    \centering
    \includegraphics[width=0.48\textwidth]{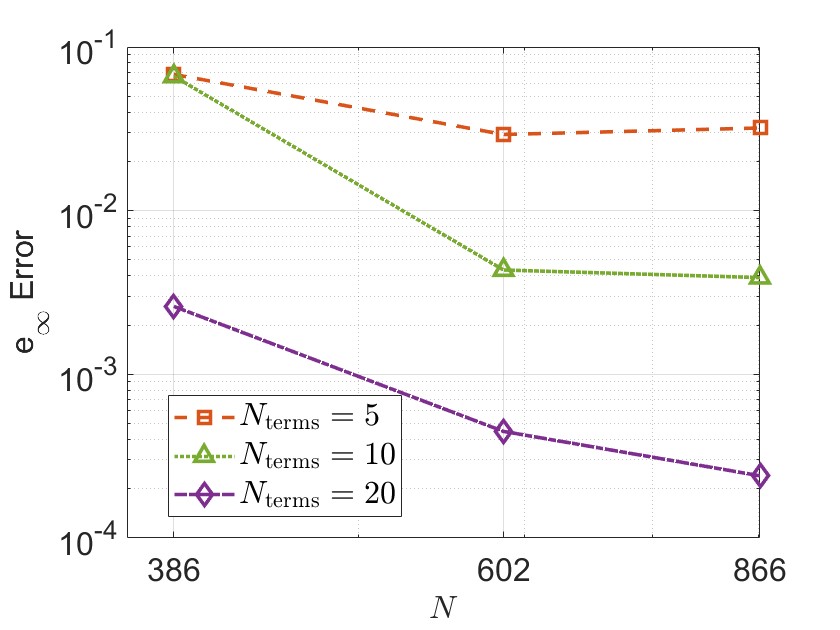}
    \hfill
    \includegraphics[width=0.48\textwidth]{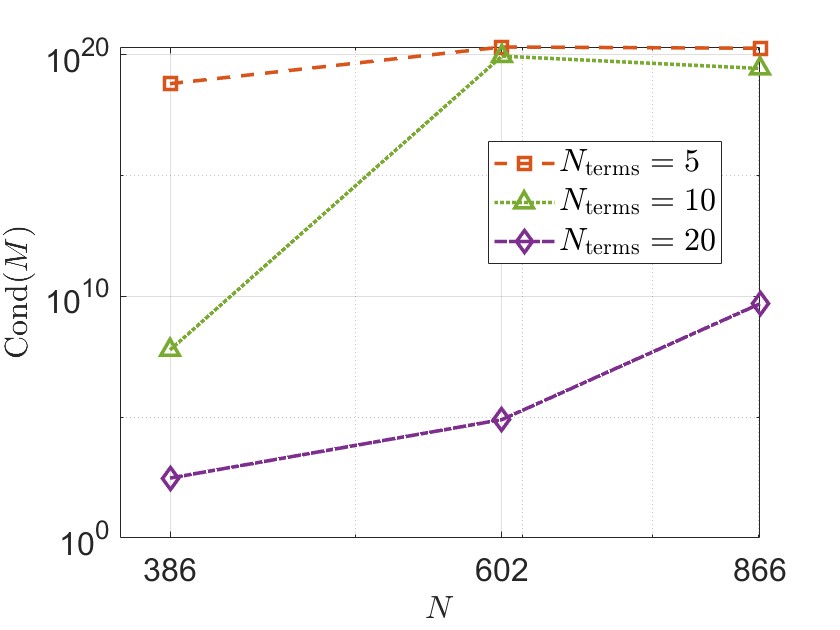}
    \caption{Dimension $d=3$. Left: $e_\infty$-error versus the number of boundary collocation points when considering different truncations of the fundamental solution up to $N_{terms}$. Right: condition number of the matrix $\mathbf{M}$.}
    \label{fig:ambas3}
\end{figure}

In Figure \ref{fig:ambas3} we show the error $e_\infty$ in terms of the number of collocation points in the boundary $N$ and the truncation $N_{terms}$ in the series of the fundamental solution. The collocation points in the boundary are computed by considering a uniform grid in each face with an increasing number of points at each edge $9,\; 11,\; 13$ which corresponds to $N=386,\; 602, \; 866$ in Figure \ref{fig:ambas3}. The error $e_\infty$ is computed as the maximum relative error on a uniform grid with $500$ points on the cube $[-5,5]^3$ that are in $\Omega$.

We observe the same qualitative behavior as in dimension $d=2$: the error decreases rapidly as both $N$ and $N_{terms}$ grow, and already for moderate values of these parameters the approximation is highly accurate. At the same time, the condition number of the matrix $\mathbf{M}$ grows rapidly with the number of points. 

\section{Proof of Theorem \ref{th_adition_theorem}}



\emph{Proof of Theorem \ref{th_adition_theorem}.} Following Arens (see \cite{Arens}) , the fundamental solution  of the Navier-Lamé operator satisfying the Kupradze condition is a $3\times 3$  symmetric  matrix $\Phi(x,y)$ satisfying $\mathrm{\Phi}(x,y)=\mathrm{\Phi}(y,x)$ and for $|x|>|y|$  given by  
\begin{equation}\label{fundamental_juan}
\mathrm{\Phi}(x,y)
=\frac{i}{\omega^2}\sum_{n=0}^\infty\sum_{m=-n}^{n}
\left(k_p\mathbf{r}_{p,n,m}(x)\overline{\mathbf{e}_{p,n,m}^t(y)}
+k_s\mathrm{R}_{s,n,m}(x)\overline{\mathrm{E}_{s,n,m}^t(y)}\right),
\end{equation}
where 
$$\mathbf{r}_{p,n,m}(x)=\nabla_x\left(h_n^{(1)}(k_p|x|)Y_n^m\left(\frac{x}{|x|}\right)\right), \hspace{0.4cm} \mathbf{e}_{p,n,m}(x)=\nabla_x\left(j_n(k_p|x|)Y_n^m\left(\frac{x}{|x|}\right)\right),$$
$$\mathrm{R}_{s,n,m}(x)=\mathcal{D}_x\left(h_n^{(1)}(k_s|x|)Y_n^m\left(\frac{x}{|x|}\right)\right), \hspace{0.4cm}
\mathrm{E}_{s,n,m}(x)=\mathcal{D}_x\left(j_n(k_s|x|)Y_n^m\left(\frac{x}{|x|}\right)\right),$$
\begin{equation*}
\mathcal{D}_x:=\left(
\begin{array}{ccc}
0&-\partial_{x_3}&\partial_{x_2}\\
\partial_{x_3}&0&-\partial_{x_1}\\
-\partial_{x_2}&\partial_{x_1}&0\\
\end{array}
\right).
\end{equation*}
The series  (\ref{fundamental_juan}) and its term by term first derivatives with respect to $|x|$ and $|y|$ are absolutely and uniformly convergent on compact subsets of $|x| > |y|$.

We define
\begin{equation}
	\label{H}
	H_{k,n,m}(x)=h_n^{(1)}(k|x|)Y_n^m\left(\frac{x}{|x|}\right), \hspace{0.3cm}  J_{k,n,m}(x)=j_n(k|x|)Y_n^m\left(\frac{x}{|x|}\right),
	\end{equation}
	
and we write
\begin{equation}
    \label{phi_matriz}
    \mathrm{\Phi}(x,y)=(\phi_{\ell j})_{\ell,j=1}^3, \; \;  \textrm{with} \;\; \phi_{\ell j}=\frac{i}{\omega^2}\sum_{n=0}^\infty\sum_{m=-n}^{n}\phi_{\ell j}^{n,m}.
\end{equation}
Then, for $\ell\neq j$ we have that
\begin{equation}
\phi_{\ell j}^{n,m}
=
k_p\
\partial_{x_\ell}H_{k_p,n,m}(x)\
\partial_{y_j}\overline{J_{k_p,n,m}(y)}
-
k_s\
\partial_{x_j}H_{k_s,n,m}(x)\
\partial_{y_\ell}\overline{J_{k_s,n,m}(y)},
\label{mixtas}
\end{equation}
whereas 
\begin{equation}\label{puras}
\phi_{\ell \ell}^{n,m}
=
k_p\
\partial_{x_\ell}H_{k_p,n,m}(x)\
\partial_{y_\ell}\overline{J_{k_p,n,m}(y)}
+k_s
\sum_{j\neq\ell}
\partial_{x_j}H_{k_s,n,m}(x)\
\partial_{y_j}\overline{J_{k_s,n,m}(y)}.
\end{equation}
In order to compute the derivatives appearing in \eqref{mixtas} and \eqref{puras} we use the spherical coordinates to write 
\begin{equation*}
x=(r\cos\varphi \sin\theta ,r\sin\varphi \sin\theta ,r\cos\theta)^t,
\; \; \; 
r>0,0\le \varphi <2\pi, 0\le \theta \le \pi.
\end{equation*}
Now we introduce the function
\begin{equation}
\label{F}
G_{k,n,m}(x):=g_n(kr)Y_n^m(\theta,\varphi)
\end{equation}
where $g_n$ could be $j_n$ or $h_n^{(1)}$.
This function satisfies the following identities, which are proved in the Appendix:
\begin{align}
\partial_{x_1}G_{k,n,m}(x)
=&
\frac{k}{2}
\left[
g_{n-1}(kr)
\left(
a_{n-1}^{m+1}Y_{n-1}^{m+1}
-b_{n-1}^{m-1}Y_{n-1}^{m-1}\right)
\right.
\nonumber
\\
&+
\left.
g_{n+1}(kr)
\left(
c_{n+1}^{m+1}Y_{n+1}^{m+1}
-d_{n+1}^{m-1}Y_{n+1}^{m-1}\right)
\right],
\label{parcial1_final} \hspace{0.25cm}
\\
\partial_{x_2}G_{k,n,m}(x)
=&
-i\frac{k}{2}
\left[
g_{n-1}(kr)
\left(
a_{n-1}^{m+1}Y_{n-1}^{m+1}
+b_{n-1}^{m-1}Y_{n-1}^{m-1}\right)
\right.
\nonumber
\\
&+
\left.
g_{n+1}(kr)
\left(
c_{n+1}^{m+1}Y_{n+1}^{m+1}
+d_{n+1}^{m-1}Y_{n+1}^{m-1}\right)
\right],
\label{parcial2_final} \hspace{0.25cm}
\\
\label{parcial3_final}
\partial_{x_3}G_{k,n,m}(x)
= &
k
\left[
g_{n-1}(kr)e_{n-1}^mY_{n-1}^m
-
g_{n+1}(kr)f_{n+1}^mY_{n+1}^m
\right],
\end{align}
where
\begin{align}
a_{n-1}^{m+1}:=&\sqrt{\frac{(n-m)(n-m-1)}{(2n+1)(2n-1)}}, \hspace{0.1cm} & b_{n-1}^{m-1}:=&\sqrt{\frac{(n+m)(n+m-1)}{(2n+1)(2n-1)}},
\label{a}
\\
c_{n+1}^{m+1}:=&\sqrt{\frac{(n+m+2)(n+m+1)}{(2n+1)(2n+3)}}, \hspace{0.1cm} & d_{n+1}^{m-1}:=&\sqrt{\frac{(n-m+2)(n-m+1)}{(2n+1)(2n+3)}}.
\label{d}
\end{align}

We start by calculating the term $\phi_{12}$ of the fundamental solution matrix $\mathrm{\Phi}(x,y)$. Unless necessary, we will omit in the notation the dependence on $\theta$ and $\varphi$.

From \eqref{mixtas} with $\ell=1$ and $j=2$ and using \eqref{parcial1_final} and \eqref{parcial2_final} with $G_{k,n,m}=H_{k_p,n,m}$ (or $J_{k_p,n,m}$) we have that
\begin{align}
& \phi_{12}^{n,m} =  
\frac{ik_p^3}{4}
\sum_{n_1\in I_n}\sum_{n_2\in I_n}
h_{n_1}^{(1)}(k_pr)j_{n_2}(k_pt)
\sum_{m_1\in I_m}\sum_{m_2\in I_m}
\tau_{12}(n_1,n_2,m_1,m_2)
Y_{n_1}^{m_1}\,\overline{Y_{n_2}^{m_2}}
\nonumber
\\
&+
\frac{ik_s^3}{4}
\sum_{n_1\in I_n}\sum_{n_2\in I_n}
h_{n_1}^{(1)}(k_sr)j_{n_2}(k_st)
\sum_{m_1\in I_m}\sum_{m_2\in I_m}
\kappa_{12}(n_1,n_2,m_1,m_2)
Y_{n_1}^{m_1}\,\overline{Y_{n_2}^{m_2}},
\label{1_elemeto_matriz} 
\end{align}
where
\begin{equation}
    \label{In}
    I_n:=\{n-1,n+1 \},
\end{equation}
$\tau_{12}(n_1,n_2,m_1,m_2)$ 
is given in Table \ref{tabla_tau}
and $\kappa_{12}(n_1,n_2,m_1,m_2)$ can be obtained from the same table as follows.
For $i,j \in \{1,2,3,4\}$ we denote by $\tau_{12}(i,j)$ the coefficient $\tau_{12}$ appearing in row $i$ and column $j$ of Table 3. With this convention, the corresponding table for the coefficients $\kappa_{12}$ can be generated from Table 3 by taking into account that $\kappa_{12}(i,j)=(-1)^{i+j}\tau_{12}(i,j)$.
\begin{table}[h]
	\centering
	\footnotesize
	\begin{tabular}{c|c|c|c|c|}
		$\tau_{12}$ &
		$\begin{array}{c}n_2=n-1\\m_2=m+1\end{array}$ &
		$\begin{array}{c}n_2=n-1\\m_2=m-1\end{array}$ &
		$\begin{array}{c}n_2=n+1\\m_2=m+1\end{array}$ &
		$\begin{array}{c}n_2=n+1\\m_2=m-1\end{array}$ \\
		\hline
		$\begin{array}{c}n_1=n-1\\m_1=m+1\end{array}$ &
		$(a_{n-1}^{m+1})^2$ &
		$a_{n-1}^{m+1}b_{n-1}^{m-1}$ &
		$a_{n-1}^{m+1}c_{n+1}^{m+1}$ &
		$a_{n-1}^{m+1}d_{n+1}^{m-1}$ \\
		\hline
		$\begin{array}{c}n_1=n-1\\m_1=m-1\end{array}$ &
		$-b_{n-1}^{m-1}a_{n-1}^{m+1}$ & 
		$-(b_{n-1}^{m-1})^2$ &
		$-b_{n-1}^{m-1}c_{n+1}^{m+1}$ &
		$-b_{n-1}^{m-1}d_{n+1}^{m-1}$ \\
		\hline
		$\begin{array}{c}n_1=n+1\\m_1=m+1\end{array}$ &
		$c_{n+1}^{m+1}a_{n-1}^{m+1}$ & 
		$c_{n+1}^{m+1}b_{n-1}^{m-1}$ &
		$(c_{n+1}^{m+1})^2$ &
		$c_{n+1}^{m+1}d_{n+1}^{m-1}$ \\
		\hline
		$\begin{array}{c}n_1=n+1\\m_1=m-1\end{array}$ &
		$-d_{n+1}^{m-1}a_{n-1}^{m+1}$ & 
		$-d_{n+1}^{m-1}b_{n-1}^{m-1}$ &
		$-d_{n+1}^{m-1}c_{n+1}^{m+1}$ &
		$-(d_{n+1}^{m-1})^2$ \\
		\hline
	\end{tabular}
	\caption{Coefficients $\tau_{12}(n_1,n_2,m_1,m_2)$.}
	\label{tabla_tau}
\end{table}

We would like to note that the order in which the factors of the products appearing in Table \ref{tabla_tau}  is written is relevant in the sense that it indicates the corresponding spherical harmonics. 
For example, the product $a_{n-1}^{m+1}d_{n+1}^{m+1}$ corresponds to harmonics $Y_{n-1}^{m+1}(\theta_x,\varphi_x)$ and $\overline{Y_{n+1}^{m+1}}(\theta_y,\varphi_y)$, 
while $c_{n+1}^{m+1}a_{n-1}^{m+1}$ corresponds to 
$Y_{n+1}^{m+1}(\theta_x,\varphi_x)$ and $\overline{Y_{n-1}^{m+1}}(\theta_y,\varphi_y)$.

Now we make the variable changes in the summation indices $n$ and $m$ appearing in \eqref{phi_matriz}, that are necessary for the spherical harmonic $Y_{n_2}^{m_2}(\theta_y,\varphi_y)$ to become $Y_{n'}^{m'}(\theta_y,\varphi_y)$. We note that $n_2\in I_n$ and $m_2\in I_m$, so four different variable changes are necessary for each index. 
More precisely,  the coefficients appearing in column 1 of Table \ref{tabla_tau}  correspond to the harmonic $Y_{n-1}^{m+1}(\theta_y,\varphi_y)$, so the variable change in this case would be $n'=n-1$ and $m'=m+1$.

After  these variable changes we group terms and write $n$ and $m$ instead of $n'$ and $m'$, we obtain
\begin{equation}
\label{phi12_final}
\phi_{12}=\frac{-1}{4\omega^2}\sum_{n=0}^\infty\sum_{m=-n}^{n}
\sum_{n_1\in K_n}\sum_{m_1\in L_m}
H_{n_1}^{-,k_p,k_s}(r,t)
\Gamma_{12}(n_1,n,m_1,m)
Y_{n_1}^{m_1}\,\overline{Y_n^m},
\end{equation}
where 
\begin{equation}
    \label{KnLm}
    K_n:=\{n-2,n,n+2 \},\qquad L_m:=\{m-2,m+2 \},
\end{equation}
$H_{n_1,n}^{-,k_p,k_s}(r,t)$ is given in (\ref{h-})
and $\Gamma_{12}(n_1,n,m_1,m)$ 
 in Table \ref{tabla_Gamma}.

\begin{table}[h]
	\centering
	\begin{tabular}{c|c|c|}
		$\Gamma_{12}$&
		$m_1=m-2$&
		$m_1=m+2$
		\\
		\hline
		$n_1=n-2$&
		$-b_{n-2}^{m-2}c_n^m$&
		$a_{n-2}^{m+2}d_n^m$
		\\
		\hline
		$n_1=n$&
		$-b_{n}^{m-2}a_n^m-d_{n}^{m-2}c_n^m$&
		$a_{n}^{m+2}b_n^m+c_{n}^{m+2}d_n^m$
		\\
		\hline 
		$n_1=n+2$&
		$-d_{n+2}^{m-2}a_n^m$&
		$c_{n+2}^{m+2}b_n^m$
		\\
		\hline
	\end{tabular}
	\caption{Coefficients $\Gamma_{12}(n_1,n,m_1,m)$.}
	\label{tabla_Gamma}
\end{table}

We would like to note that $H_{n_1,n}^{-,k_p,k_s}=0$ when $k_p=k_s$.

Arguing as we have done with $\phi_{12}$ to obtain \eqref{phi12_final}, similar expressions for $\phi_{13}$ and  $\phi_{23}$ can be derived. More specifically
\begin{align}
\phi_{13}= &  \frac{i}{2\omega^2}\sum_{n=0}^\infty\sum_{m=-n}^{n}
\sum_{n_1\in K_n}\sum_{m_1\in I_m}
H_{n_1,n}^{-,k_p,k_s}(r,t)
\Gamma_{13}(n_1,n,m_1,m)
Y_{n_1}^{m_1}\,\overline{Y_n^m},
\label{phi13_final}
\\
\label{phi23_final}
\phi_{23}= & \frac{1}{2\omega^2}\sum_{n=0}^\infty\sum_{m=-n}^{n}
\sum_{n_1\in K_n}\sum_{m_1\in I_m}
H_{n_1,n}^{-,k_p,k_s}(r,t)
\Gamma_{23}(n_1,n,m_1,m)
Y_{n_1}^{m_1}\,\overline{Y_n^m},
\end{align}

where $K_n$ is given in \eqref{KnLm}, $I_m$ in \eqref{In}, $H_{n_1,n}^{-,k_p,k_s}$ in \eqref{h-}, $\Gamma_{13}(n_1,n,m_1,m)$ in Table \ref{tabla_Gamma13} and $\Gamma_{23}(n_1,n,m_1,m)$ in Table \ref{tabla_Gamma23}, with
\begin{equation}\label{e}
e_{n-1}^m := \sqrt{\frac{(n+m)(n-m)}{(2n+1)(2n-1)}},\hspace{0.5cm}  f_{n+1}^m := \sqrt{\frac{(n+m+1)(n-m+1)}{(2n+1)(2n+3)}}.
\end{equation}

\begin{table}[h]
	\centering
	\begin{tabular}{c|c|c|}
		$\Gamma_{13}$&
		$m_1=m-1$&
		$m_1=m+1$
		\\
		\hline
		$n_1=n-2$&
		$b_{n-2}^{m-1}f_n^m$&
		$-a_{n-2}^{m+1}f_n^m$
		\\
		\hline
		$n_1=n$&
		$-b_{n}^{m-1}e_n^m+d_{n}^{m-1}f_n^m$&
		$a_{n}^{m+1}e_n^m-c_{n}^{m+1}f_n^m$
		\\
		\hline 
		$n_1=n+2$&
		$-d_{n+2}^{m-1}e_n^m$&
		$c_{n+2}^{m+1}e_n^m$
		\\
		\hline
	\end{tabular}
	\caption{Coefficients $\Gamma_{13}(n_1,n,m_1,m)$.}
	\label{tabla_Gamma13}
\end{table}

\begin{table}[h]
	\centering
	\begin{tabular}{c|c|c|}
		$\Gamma_{23}$&
		$m_1=m-1$&
		$m_1=m+1$
		\\
		\hline
		$n_1=n-2$&
		$-b_{n-2}^{m-1}f_n^m$&
		$-a_{n-2}^{m+1}f_n^m$
		\\
		\hline
		$n_1=n$&
		$b_{n}^{m-1}e_n^m-d_{n}^{m-1}f_n^m$&
		$a_{n}^{m+1}e_n^m-c_{n}^{m+1}f_n^m$
		\\
		\hline 
		$n_1=n+2$&
		$d_{n+2}^{m-1}e_n^m$&
		$c_{n+2}^{m+1}e_n^m$
		\\
		\hline
	\end{tabular}
	\caption{Coefficients $\Gamma_{23}(n_1,n,m_1,m)$.}
	\label{tabla_Gamma23}
\end{table}

We note that the spherical harmonics $Y_{n_1}^{m_1}$ appearing in \eqref{phi12_final}, \eqref{phi13_final} and \eqref{phi23_final} are defined only for $n_1\ge 0$ and $-n_1\le m_1\le n_1$.
This poses no difficulty, since one can verify that the factors $\Gamma_{12}(n_1,n,m_1,m),\ \Gamma_{13}(n_1,n,m_1,m)$ and $\Gamma_{23}(n_1,n,m_1,m)$ accompanying the spherical harmonics vanish outside this range.

As $\mathrm{\Phi}$ is a symmetric matrix, we only need to find an expression for $\phi_{\ell\ell}$, $\ell=1,2,3$. 

Let's start with $\phi_{11}$. 
From \eqref{phi_matriz} with $\ell=j=1$ we have that
\begin{equation}
\label{phi11}
\phi_{11}=\frac{i}{\omega^2}\sum_{n=0}^\infty\sum_{m=-n}^{n}\phi_{11}^{n,m},
\end{equation}
where $\phi_{11}^{n,m}$ is given by \eqref{puras}.
Using \eqref{parcial1_final}, \eqref{parcial2_final} and \eqref{parcial3_final} with $G_{k,n,m}=H_{k_p,n,m}$ (or $J_{k_p,n,m}$) we have that
\begin{align}
& \phi_{11}^{n,m}= \nonumber
 \frac{k_p^3}{4}
\sum_{n_1\in I_n}\sum_{n_2\in I_n}
h_{n_1}^{(1)}(k_pr)j_{n_2}(k_pt)
\sum_{m_1\in I_m}\sum_{m_2\in I_m}
\tau(n_1,n_2,m_1,m_2)
Y_{n_1}^{m_1}\,\overline{Y_{n_2}^{m_2}} \nonumber
\\
&+
\frac{k_s^3}{4}
\sum_{n_1\in I_n}\sum_{n_2\in I_n}
h_{n_1}^{(1)}(k_sr)j_{n_2}(k_st)
\sum_{m_1\in I_m}\sum_{m_2\in I_m}
\kappa(n_1,n_2,m_1,m_2)
Y_{n_1}^{m_1}\,\overline{Y_{n_2}^{m_2}}  \nonumber
\\
&+
k_s^3
\sum_{n_1\in I_n}\sum_{n_2\in I_n}
h_{n_1}^{(1)}(k_sr)j_{n_2}(k_st)
\chi(n_1,n_2,m,m)
Y_{n_1}^{m}\,\overline{Y_{n_2}^{m}},  \label{phi_11_juan}
\end{align}
where $I_n$ is given in \eqref{In}; the coefficients $\tau(n_1,n_2,m_1,m_2)$ and $\chi(n_1,n_2,m_1,m_2)$ are listed in Tables \ref{tabla_tau_11} and \ref{tabla_chi11}, respectively; and $\kappa(n_1,n_2,m_1,m_2)$ can be obtained from Table \ref{tabla_tau_11} using the convention introduced to derive the coefficients $\kappa_{12}$ in \eqref{1_elemeto_matriz}, together with the relation $\kappa(i,j)=(-1)^{i+j}\tau(i,j)$.
\begin{table}[h]
	\centering
	\footnotesize 
	\begin{tabular}{c|c|c|c|c|}
		$\tau$ &
		$\begin{array}{c}n_2=n-1\\m_2=m+1\end{array}$ &
		$\begin{array}{c}n_2=n-1\\m_2=m-1\end{array}$ &
		$\begin{array}{c}n_2=n+1\\m_2=m+1\end{array}$ &
		$\begin{array}{c}n_2=n+1\\m_2=m-1\end{array}$ \\
		\hline
		$\begin{array}{c}n_1=n-1\\m_1=m+1\end{array}$ &
		$(a_{n-1}^{m+1})^2$ &
		$-a_{n-1}^{m+1}b_{n-1}^{m-1}$ &
		$a_{n-1}^{m+1}c_{n+1}^{m+1}$ &
		$-a_{n-1}^{m+1}d_{n+1}^{m-1}$ \\
		\hline
		$\begin{array}{c}n_1=n-1\\m_1=m-1\end{array}$ &
		$-b_{n-1}^{m-1}a_{n-1}^{m+1}$ &
		$(b_{n-1}^{m-1})^2$ &
		$-b_{n-1}^{m-1}c_{n+1}^{m+1}$ &
		$b_{n-1}^{m-1}d_{n+1}^{m-1}$ \\
		\hline
		$\begin{array}{c}n_1=n+1\\m_1=m+1\end{array}$ &
		$c_{n+1}^{m+1}a_{n-1}^{m+1}$ &
		$-c_{n+1}^{m+1}b_{n-1}^{m-1}$ &
		$(c_{n+1}^{m+1})^2$ &
		$-c_{n+1}^{m+1}d_{n+1}^{m-1}$ \\
		\hline
		$\begin{array}{c}n_1=n+1\\m_1=m-1\end{array}$ &
		$-d_{n+1}^{m-1}a_{n-1}^{m+1}$ &
		$d_{n+1}^{m-1}b_{n-1}^{m-1}$ &
		$-d_{n+1}^{m-1}c_{n+1}^{m+1}$ &
		$(d_{n+1}^{m-1})^2$ \\
		\hline
	\end{tabular}
	\caption{Coefficients $\tau_{11}(n_1,n_2,m_1,m_2)$.}
	\label{tabla_tau_11}
\end{table}
\begin{table}[h]
	\centering
	\begin{tabular}{c|c|c|}
		$\chi$&
		$n_2=n-1$&
		$n_2=n+1$
		\\
		\hline
		$n_1=n-1$&
		$(e_{n-1}^m)^2$&
		$-e_{n-1}^mf_{n+1}^m$
		\\
		\hline
		$n_1=n+1$&
		$-f_{n+1}^me_{n-1}^m$&
		$(f_{n+1}^m)^2$
		\\
		\hline
	\end{tabular}
	\caption{Coefficients $\chi_{11}(n_1,n_2,m,m)$.}
	\label{tabla_chi11}
\end{table}

After making the necessary variable changes in the summation indices $n$ and $m$ in \eqref{phi_11_juan} so that the spherical harmonics $Y_{n}^{m}(\theta_y,\varphi_y)$ always appear, we have
\begin{align}
& \phi_{11}=\frac{i}{\omega^2}\sum_{n=0}^\infty\sum_{m=-n}^{n}
H_{n,n}^{k_s}(r,t)
Y_{n}^{m}\,\overline{Y_n^m}
\nonumber
\\
&+\frac{i}{4\omega^2}\sum_{n=0}^\infty\sum_{m=-n}^{n}
\sum_{n_1\in K_n}\sum_{m_1\in K_m}
H_{n_1,n}^{-,k_p,k_s}(r,t)
\Gamma_{11}(n_1,n,m_1,m)
Y_{n_1}^{m_1}\,\overline{Y_n^m},
\label{phi11_final}
\end{align}
where $H_{n,n}^{k_s}$ is given in \eqref{h0},
$K_n$  in \eqref{KnLm}, $H_{n_1,n}^{-,k_p,k_s}$ in \eqref{h-} and $\Gamma_{11}(n_1,n,m_1,m)$ in Table \ref{tabla_Gamma11}.

\begin{table}[h]
	\centering
  \resizebox{\textwidth}{!}{%
	\begin{tabular}{c|c|c|c|}
		$\Gamma_{11}$&
		$m_1=m-2$&
		$m_1=m$&
		$m_1=m+2$
		\\
		\hline
		$n_1=n-2$&
		$-b_{n-2}^{m-2}c_n^m$&
		$2a_{n-2}^{m}c_n^m$&
		$-a_{n-2}^{m+2}d_n^m$
		\\
		\hline
		$n_1=n$&
		$-b_{n}^{m-2}a_n^m-d_{n}^{m-2}c_n^m$&
		$(a_{n}^m)^2+(b_{n}^m)^2+(c_{n}^m)^2+(d_{n}^m)^2$&
		$-a_{n}^{m+2}b_n^m-c_{n}^{m+2}d_n^m$
		\\
		\hline 
		$n_1=n+2$&
		$-d_{n+2}^{m-2}a_n^m$&
		$2c_{n+2}^{m}a_n^m$&
		$-c_{n+2}^{m+2}b_n^m$
		\\
		\hline
	\end{tabular}
 }
	\caption{Coefficients $\Gamma_{11}(n_1,n,m_1,m)$.}
	\label{tabla_Gamma11}
\end{table}

Arguing in a similar way we get
\begin{equation*}
\label{phi22}
\phi_{22}=\frac{i}{\omega^2}\sum_{n=0}^\infty\sum_{m=-n}^{n}\phi_{22}^{n,m},
\end{equation*}
with
\begin{align*}
\phi_{22}^{n,m}=&
\frac{k_p^3}{4}
\sum_{n_1\in I_n}\sum_{n_2\in I_n}
h_{n_1}^{(1)}(k_pr)j_{n_2}(k_pt)
\sum_{m_1\in I_m}\sum_{m_2\in I_m}
\kappa_{11}(n_1,n_2,m_1,m_2)
Y_{n_1}^{m_1}\,\overline{Y_{n_2}^{m_2}}
\\
&+
\frac{k_s^3}{4}
\sum_{n_1\in I_n}\sum_{n_2\in I_n}
h_{n_1}^{(1)}(k_sr)j_{n_2}(k_st)
\sum_{m_1\in I_m}\sum_{m_2\in I_m}
\tau_{11}(n_1,n_2,m_1,m_2)
Y_{n_1}^{m_1}\,\overline{Y_{n_2}^{m_2}}
\\
&+
k_s^3
\sum_{n_1\in I_n}\sum_{n_2\in I_n}
h_{n_1}^{(1)}(k_sr)j_{n_2}(k_st)
\chi(n_1,n_2,m,m)
Y_{n_1}^{m}\,\overline{Y_{n_2}^{m}}.
\end{align*} 
where $I_n$ is given in \eqref{In}. 

From here, we obtain
\begin{align}
& \phi_{22}=\frac{i}{\omega^2}\sum_{n=0}^\infty\sum_{m=-n}^{n}
H_{n,n}^{k_s}(r,t)
Y_{n}^{m}\,\overline{Y_n^m}
\nonumber
\\
&+\frac{i}{4\omega^2}\sum_{n=0}^\infty\sum_{m=-n}^{n}
\sum_{n_1\in K_n}\sum_{m_1\in K_m}
H_{n_1,n}^{-,k_p,k_s}(r,t)
\Gamma_{22}(n_1,n,m_1,m)
Y_{n_1}^{m_1}\,\overline{Y_n^m},
\label{phi22_final}
\end{align}
with $K_n$ given in \eqref{KnLm}, $H_{n,n}^{k_s}$ in \eqref{h0}, $H_{n_1,n}^{-,k_p,k_s}$ in \eqref{h-} and,
\begin{equation*}
\Gamma_{22}(n_1,n,m_1,m)=
\left\{
\begin{array}{rll}
\Gamma_{11}(n_1,n,m_1,m)&& \text{if }m_1=m,
\\
-\Gamma_{11}(n_1,n,m_1,m)&& \text{if }m_1\neq m.
\end{array}
\right.
\end{equation*}

Analogously,
\begin{equation*}
\label{phi33}
\phi_{33}=\frac{i}{\omega^2}\sum_{n=0}^\infty\sum_{m=-n}^{n}\phi_{33}^{n,m},
\end{equation*}
with
\begin{align*}
\phi_{33}^{n,m}=&
k_p^3
\sum_{n_1\in I_n}\sum_{n_2\in I_n}
h_{n_1}(k_pr)j_{n_2}(k_pt)
\chi_{11}(n_1,n_2,m,m)
Y_{n_1}^{m}\,\overline{Y_{n_2}^{m}}
\\
&+
\frac{k_s^3}{4}
\sum_{n_1\in I_n}\sum_{n_2\in I_n}
h_{n_1}(k_sr)j_{n_2}(k_st)
\sum_{m_1\in I_m}\sum_{m_2\in I_m}
\tau_{11}(n_1,n_2,m_1,m_2)
Y_{n_1}^{m_1}\,\overline{Y_{n_2}^{m_2}}
\\
&+
\frac{k_s^3}{4}
\sum_{n_1\in I_n}\sum_{n_2\in I_n}
h_{n_1}(k_sr)j_{n_2}(k_st)
\sum_{m_1\in I_m}\sum_{m_2\in I_m}
\kappa_{11}(n_1,n_2,m_1,m_2)
Y_{n_1}^{m_1}\,\overline{Y_{n_2}^{m_2}},
\end{align*}

where $I_n$ is given in \eqref{In}. 

After making the necessary variable changes and grouping terms, all terms cancel out except those associated with the spherical harmonics $Y_{n_1}^{m}(\theta_x,\varphi_x)$ with $n_1\in K_n$.
Therefore, we get
\begin{align}
& \phi_{33}=\frac{i}{\omega^2}\sum_{n=0}^\infty\sum_{m=-n}^{n}
H_{n,n}^{k_s}(r,t)
Y_{n}^{m}\,\overline{Y_n^m}
\nonumber
\\
&+\frac{i}{\omega^2}\sum_{n=0}^\infty\sum_{m=-n}^{n}
\sum_{n_1\in K_n}
H_{n_1,n}^{-,k_p,k_s}(r,t)
\Gamma_{33}(n_1,n,m,m)
Y_{n_1}^{m}\,\overline{Y_n^m},
\label{phi33_final}
\end{align}
with $K_n$ given in \eqref{KnLm}, $H_{n,n}^{k_s}$ in \eqref{h0}, $H_{n_1,n}^{-,k_p,k_s}$ in \eqref{h-} and,
\begin{equation*}
\Gamma_{33}(n_1,n,m,m)=
\left\{
\begin{array}{lll}
-e_{n-2}^mf_n^m&& \text{if }n_1=n-2,
\\
(e_{n}^m)^2+(f_{n}^m)^2&& \text{if }n_1=n,
\\
-f_{n+2}^me_n^m&& \text{if }n_1=n+2.
\end{array}
\right.
\end{equation*}

As before, we note that the fact that the spherical harmonics $Y_{n_1}^{m_1}$ appearing in \eqref{phi11_final}, \eqref{phi22_final}, and \eqref{phi33_final} may involve indices $n_1$ and $m_1$ outside their range of definition does not pose any difficulty, since one can verify that the factors $\Gamma_{11}(n_1,n,m_1,m)$, $\Gamma_{22}(n_1,n,m_1,m)$, and $\Gamma_{33}(n_1,n,m_1,m)$ accompanying the spherical harmonics vanish in those cases.

Using \eqref{phi12_final}, \eqref{phi13_final}, \eqref{phi23_final}, \eqref{phi11_final}, \eqref{phi22_final} and \eqref{phi33_final}, the result follows from \eqref{phi_matriz}.

\hfill$\square$

\section{Appendix}
In this appendix we will prove the identities \eqref{parcial1_final}, \eqref{parcial2_final} and \eqref{parcial3_final}.

We consider the following normalization for the spherical harmonics:
\begin{equation}
\label{gamma_n^m}
Y_n^m(\theta,\varphi):=\gamma_n^me^{im\varphi}P_n^{m}(\cos\theta), \hspace{0.5cm}\gamma_n^m:=\sqrt{\frac{(2n+1)(n-m)!}{4\pi(n+m)!}}.
\end{equation}
Here $P_n^{m}$ is the associated Legendre functions of degree $n$ and order $m$ defined by
\begin{equation}
\label{legendre}
P_n^{m}(z):=(-1)^m(1-z^2)^\frac{m}{2}\frac{d^mP_n(z)}{dz^m}, 
\qquad
z\in[-1,1],\ 0\le m\le n,
\end{equation}
where $P_n$ is the Legendre polynomial of degree $n$.
We use the standard convention that for $m>0$,
\begin{equation}
\label{Pnm_negativo}
P_n^{-m}(z):=(-1)^m\frac{(n-m)!}{(n+m)!}\,P_n^m(z).
\end{equation}

With this convention, the identity (see \cite[pp. 334, 8.5.4]{AS1967})
\begin{equation*}
\label{derivadaP}
(z^2-1)\frac{dP_n^m(z)}{dz}=nzP_n^m(z)-(n+m)P_{n-1}^m(z), \qquad z\in[-1,1],
\end{equation*}
also holds for negative orders, so it holds for $|m|\le n$.
We would like to note that this expression is valid for $|m|=n$, since from \eqref{legendre} we have that $P_{n-1}^n=0$.

Using this identity, we have that

\begin{equation}\label{parcial_theta}
\partial_\theta Y_n^m
=\gamma_n^m e^{im\varphi}\frac{n\cos\theta P_n^m(\cos\theta)-(n+m)P_{n-1}^m(\cos\theta)}{\sin\theta}.
\end{equation}

We will need the following identities  involving the associated Legendre functions, which can be obtained from the usual relations of these functions, see e.g. \cite[chapter 8]{AS1967}, \cite[chapter 7]{Lebedev} or \cite[chapter 2]{Edmons}:
\begin{equation}
    (n-m)\frac{P_n^m(\cos\theta)}{\sin\theta}-(n+m)\cot\theta P_{n-1}^m(\cos\theta)
=P_{n-1}^{m+1}(\cos\theta), \label{juan_1}
\end{equation}
\begin{equation}
    (n+m)\frac{P_n^m(\cos\theta)}{\sin\theta}-(n+m)\cot\theta P_{n-1}^m (\cos \theta)
=-(n+m)(n+m-1)P_{n-1}^{m-1}(\cos\theta), \label{juan_2}
\end{equation}
\begin{equation}
\begin{split}
    (2n+1)\sin\theta P_n^m(\cos\theta)
    &=P_{n-1}^{m+1}(\cos\theta)-P_{n+1}^{m+1}(\cos\theta)
    \\
    &=(n-m+1)(n-m+2)P_{n+1}^{m-1}(\cos\theta)
    \\
    &\hspace{0.5cm}-(n+m-1)(n+m)P_{n-1}^{m-1}(\cos\theta), 
\end{split}
\label{juan_3}
\end{equation}
\begin{equation}
  \label{juan_4}
(2n+1)\cos\theta P_n^m(\cos\theta)
=(n-m+1)P_{n+1}^m(\cos\theta)
+(n+m)P_{n-1}^m(\cos\theta).  
\end{equation}

On the other hand, consider the functions $G_{k,n,m}$ and $g_n$ given in \eqref{F}.
Taking into account that (see \cite{AS1967}),
$$ \frac{g_n(r)}{r}=\frac{g_{n-1}(r)+g_{n+1}(r)}{2n+1}, \hspace{1cm}
g_n'(r)=\frac{ng_{n-1}(r)-(n+1)g_{n+1}(r)}{2n+1}, $$
we obtain, for $j=1,2,3$,
\begin{equation}
\label{parcialj_ab}
\partial_{x_j}G_{k,n,m}(x)
=
\frac{k}{2n+1}
\left[
g_{n-1}(kr)\alpha_{j,n,m}(\theta,\varphi)
+
g_{n+1}(kr)\beta_{j,n,m}(\theta,\varphi)
\right],
\end{equation}
where
\begin{align}
\label{alpha_jnm}
\alpha_{j,n,m}(\theta,\varphi)
=&nY_n^m\partial_{x_j}r
+r\partial_\theta Y_n^m\partial_{x_j}\theta
+r\partial_\varphi Y_n^m\partial_{x_j}\varphi,
\\
\label{beta_jnm}
\beta_{j,n,m}(\theta,\varphi)
=&-(n+1)Y_n^m\partial_{x_j}r
+r\partial_\theta Y_n^m\partial_{x_j}\theta
+r\partial_\varphi Y_n^m\partial_{x_j}\varphi.
\end{align}
For convenience, we omit the dependence of the involved functions on $\theta$ and $\varphi$.

Writing $\cos \varphi = \frac{e^{i \varphi } +e^{-i \varphi }  }{2} $ and $\sin \varphi = \frac{-i(e^{i \varphi } -e^{-i \varphi })  }{2} $ and using \eqref{gamma_n^m} and  \eqref{parcial_theta} gives
\begin{align*}
\alpha_{1,n,m} & =\frac{\gamma_n^m}{2} e^{i(m+1)}\left( (n-m)\frac{P_n^m(\cos\theta)}{\sin\theta}-(n+m)\cot\theta P_{n-1}^m(\cos\theta)  \right) 
\\
& +\frac{\gamma_n^m}{2} e^{i(m-1)}\left( (n+m)\frac{P_n^m(\cos\theta)}{\sin\theta}-(n+m)\cot\theta P_{n-1}^m(\cos\theta) \right).
\end{align*}
From \eqref{juan_1}, \eqref{juan_2} and \eqref{gamma_n^m}    we get
\begin{equation}
\label{alpha1nmf}
\alpha_{1,n,m}
=\frac{\gamma_n^m}{2\gamma_{n-1}^{m+1}}Y_{n-1}^{m+1}
-\frac{\gamma_n^m(n+m)(n+m-1)}{2\gamma_{n-1}^{m-1}}Y_{n-1}^{m-1}.
\end{equation}

In a similar way we obtain
$$\beta_{1,n,m} =\alpha_{1,n,m} - \frac{\gamma_n^m}{2} (2n+1) \sin \theta P_n^m(\cos\theta)  \left(e^{i(m+1)\varphi}+ e^{i(m-1)\varphi}  \right). $$
From \eqref{juan_3} and \eqref{gamma_n^m}    we get
\begin{equation}
\label{beta1nmf}
\beta_{1,n,m}
=\frac{\gamma_n^m}{2\gamma_{n+1}^{m+1}}Y_{n+1}^{m+1}
-\frac{\gamma_n^m(n-m+1)(n-m+2)}{2\gamma_{n+1}^{m-1}}Y_{n+1}^{m-1}.
\end{equation}
Finally, inserting \eqref{alpha1nmf} and \eqref{beta1nmf} in \eqref{parcialj_ab}, and using \eqref{gamma_n^m}, we obtain \eqref{parcial1_final}.

In a similar way we obtain \eqref{parcial2_final}.

By using \eqref{gamma_n^m}, \eqref{parcial_theta} and \eqref{juan_4} we have
\begin{equation*}
\alpha_{3,n,m}
=\frac{\gamma_n^m(n+m)}{\gamma_{n-1}^m}Y_{n-1}^m,
\qquad
\beta_{3,n,m}
=-\frac{\gamma_n^m(n-m+1)}{\gamma_{n+1}^m}Y_{n+1}^m.
\end{equation*}

Finally, using these identities in \eqref{parcialj_ab} we obtain   \eqref{parcial1_final}-\eqref{parcial3_final}.

\section{Conclusions}

We have derived a new addition formula for the fundamental solution of the three-dimensional Navier–Lamé system satisfying the Kupradze radiation condition at infinity. This result generalizes the two-dimensional formulation presented in \cite{BRVW2024}. As an application, we demonstrate how these formulas can be employed within the MFS to solve exterior boundary value problems, yielding accurate approximations when the series expansion of the fundamental solution is truncated with a few terms.  It would be interesting to find convergence formulas, at least in simple geometries, such as those known for the Laplace and Hemholtz equation in a ball (see \cite{Bogomolny1985}, \cite{Katsurada1990}, and \cite{BarnettBetcke2008}). This is the subject of ongoing research.





\Addresses
\end{document}